\documentstyle[12pt]{article}

\begin{document}

\newtheorem{theorem}{Theorem}[section]
\newtheorem{corollary}{Corrolary}[section]
\newtheorem{lemma}{Lemma}[section]
\newtheorem{conjecture}{Conjecture}[section]
\newtheorem{convention}{Convention}[section]
\newtheorem{remark}{Remark}[section]
\newtheorem{definition}{Definition}[section]
\newtheorem{proposition}{Proposition}[section]
\renewcommand{\baselinestretch}{1.4}
\newcommand{\binom}{ { { {\rm {wt}}  a} \choose i } }
\newcommand{\wta}{{\rm {wt} }  a }
\newcommand{\R}{\mbox {Res} _{z} }
\newcommand{\wtb}{{\rm {wt} }  b }
\newcommand{\bea}{\begin{eqnarray}}
\newcommand{\eea}{\end{eqnarray}}
\newcommand{\nn}{\nonumber \\}
\newcommand{\be}{\begin {equation}}
\newcommand{\ee}{\end{equation}}
\newcommand{\g}{\bf g}
\newcommand{\hg}{\hat {\bf g} }
\newcommand{\h}{\bf h}
\newcommand{\hh}{\hat {\bf h} }
\newcommand{\n}{\bf n}
\newcommand{\hn} {\hat {\bf n}}
\baselineskip=12pt

\title{ {Tensor Product of Vertex Operator Algebras } }
\date{17.8.1995}
\author{Antun Milas}
\maketitle

\begin{abstract}
Let $V$=$V_1 \otimes V_2$ be a tensor product of VOAs. Using Zhu
theory we discuss the theory of representations of $V$ (associative
algebra, modules and 'fusion rules'). We prove
that this theory is more or less the same as representation theory
of tensor product of the associative algebras.
\end{abstract}

\pagenumbering{arabic}

\section{Introduction}
In this paper we discuss some things which are closely related with so
called aksiomatic approach to VOA [FHL], [FLM], [B], [Li]
[DL]. There are lot of literature in theory
of VOAs, representation theory etc. It seems that Borchards [B] first make a structure of
VOA on tensor product on very natural way. It's not worth to mention that
construction of tensor product of modules is rather complicated [Li]
and use the theory of intertwining operators. 

In this paper we
gave some general results about tensor product of VOAs. We find that
associative algebra ([Z], [FZ]) have a good functorial property related
to tensor
product of two VOAs (Theorem \ref {izo}). 

The other aspect which is considered is representation theory (classification of modules
and calculating the 'fusion rules'). Using the upper property of associative
algebra we prove the classification theorem (see [FHL])  and discuss
whether tensor product of two VOAs is a rational VOA.

Using the structure of bimodules (section 5.) we present the proof 
of fusion rules formula for $V=V_1 \otimes V_2$.

It's easy to make the similar theory in case of tensor product of
two superalgebras [KWn], [Li1].
Also the similar results are valid in the classification theory
of twisted modules of tensor product of VOA [DLM].

The next step closely related with this theory would be 
consideration of {\em dual pairs} in VOAs (see [FZ]), or {\em coset}
models (WZNW models, minimal models). 
If we have VOA $V$ and $V_1$ and $V_2$ subalgebras which 'commute' then we 
said that $(V_1,V_2)$ form a dual pair for $V$--module $M$ if 
$M= \oplus_i M_i \otimes N_i,$ where $M_i$($N_i$) are $V_1$ ($V_2$)--modules.
We are working on this problematic.

\section{Preliminaries}

Vertex operator algebras is mathematical counterpart of Conformal field 
theory (CFT). Today there are various definitions of VOA and representations
of VOA according to graduation, Virasoro element etc..So we're going to
define VOA and modules in generality what we need.

\subsection{Vertex operator algebras and modules }

\begin{definition} \label{def.v}
A {\em vertex operator algebra} is a
${\bf Z}$--graded vector space $V=\bigoplus_
{n \in{\bf Z}}V_n$ with
a sequence of linear operators
$\{ a(n) \mid n\in  {\bf Z} \}\subset End\ V$ associated to every $a\in V$,
such that for fixed
$a,b\in V$, $a(n)b=0$ for $n$ sufficiently large.
We call the
 generating series
$Y(a,z)=\sum_{n\in {\bf Z}}a(n)z^{-n-1}\in (End\ V)[[z,z^{-1}]]$,
 {\em vertex operators} associated to $a$,
satisfy the following axioms:
\item[{\bf (V1)}]  $Y(a,z)=0$ iff $a=0$.

\item[{\bf (V2)}] There is a {\em vacuum} vector, which we denote by $\bf1$,
such that
$$Y({\bf 1},z)=I_V\,(I_V\, \mbox{is the identity of}\,\, End\,V). $$

\item[{\bf (V3)}]
There is a special element $\omega \in V$
(called the  {\em Virasoro\ element}), whose
 vertex operator we write in the form
$$Y(\omega,z)=\sum_{n\in {\bf Z}}\omega(n)z^{-n-1}=
\sum_ {n\in {\bf Z}}L_n z^{-n-2},$$ such that
$$L_0\mid_{V_n}=nI\mid_{V_n},$$
\bea
Y(L_{-1}a,z)=\frac{d}{dz}Y(a,z)\ \ \mbox{for every } a \in V ,
\eea
\bea
[L_m,L_n]=(m-n)L_{m+n}+\delta_{m+n,0}\frac{m^3-m}{12}c, \label{eq_vir}
\eea
where $c$ is some constant in ${\bf{C}}$,
which is called the {\em rank} of $V$.

\item[{\bf (V4)}] The {\em Jacobi identity} holds, i.e.
\bea
& &z_{0}^{-1}\delta\left(\frac{z_{1}-z_{2}}{z_{0}}\right)Y(a,z_{1})Y(b,z_{2})
-z_{0}^{-1}\delta\left(\frac{-z_{2}+z_{1}}{z_{0}}\right)Y(b,z_{2})Y(a,z_{1})
\nonumber\\
&=&z_{2}^{-1}\delta\left(\frac{z_{1}-z_{0}}{z_{2}}\right)Y(Y(a,z_{0})b,z_{2})
\eea
for any $a, b \in V$.
\end{definition}

\begin{definition}
Let $V_1$ and $V_2$ be a two VOAs then mapping
$ f \ : \ V_1 \rightarrow V_2$
is homomorphism of VOAs if F is a homomorphism of vector spaces
and if $f(Y(a,z)b)=Y(f(a),z)f(b)$ for every $a,b \in V$. $f$ is a
isomorphism of VOAs if $f$ is bijection.
\end{definition}

Let $V_1$ and $V_2$ be a two VOAs. Then the space $V_1 \oplus V_2$
has a natural VOA structure.
The subspace $I$ of $V$ is called ideal if
$
Y(a,z)b \in I[[z,z^{-1}]]$
for every $a \in V, b \in I$.
Given an ideal $I$ in $V$ such that ${\bf 1} \notin I$, $\omega \notin I$,
the quotient $V/I$ admits a natural VOA structure( see [FZ] ).
VOA is simple if doesn't have nontrivial ideals. Also, we define semisimple
VOA as VOA which is direct sum of simple VOAs.

\begin{definition}
Given an VOA $V$, {\em a representation of $V$
(or V-module)} is a ${\bf Z}_{+} $-graded vector space
 $M=\bigoplus_{n \in {\bf Z}_{+}}M_{n}$ and
a linear map $$V\longrightarrow (End\ M)[[z,z^{-1}]],$$
$$a\longmapsto Y_M (a,z)=\sum_{n\in {\bf Z}} a(n)z^{-n-1},$$
satisfying
\item[{\bf (M1)}]  $a(n)M_m \subset M_{m+\deg a -n-1}\,\,\, \mbox{for every
homogeneous element a.}$
\item [{\bf (M2)}]  $Y_M(1,z)=I_M,$ and setting
$Y_M(\omega,z)=\sum_{n\in {\bf Z} }L_nz^{-n-2},$ we have
\bea
&&[L_m,L_n]=(m-n)L_{m+n}+\delta_{m+n,0}\frac{m^3-m}{12}c,\ nonumber \\
&&Y_M(L_{-1}a,z)=\frac{d}{dz}Y_M(a,z)\,\, \nonumber
\eea
for every $\,\,\,a \in V.$

\item[{\bf (M3)}]  The {\em Jacobi identity} holds, i.e.
\bea
&&z_{0}^{-1}\delta\left(\frac{z_{1}-z_{2}}{z_{0}}\right)
Y_M (a,z_{1})Y_M (b,z_{2})
-z_{0}^{-1}\delta\left(\frac{-z_{2}+z_{1}}{z_{0}}\right)
Y_M(b,z_{2})Y_M(a,z_{1})
\nonumber \\
&&z_{2}^{-1}\delta\left(\frac{z_{1}-z_{0}}{z_{2}}\right)
Y_M(Y(a,z_{0})b,z_{2})
\eea
for any $a, b \in V$.
\label{definition_module}

\end{definition}

\begin{remark}
In [FHL]  the definition of $V$--module is 
stronger, with the following additional aksiom:
\item[{\bf(M0)}]
Dimension of $M_0(n)$ (graded subspace) is finite for every $n$.

\end{remark}

\subsection{Intertwinig operators}

\begin{definition}

Let $M^1,M^2$ and $M^3$ be three $V$--modules. An intertwinig operator of
type   $\left(\begin{array}{c}M^3 \\ M^1, \ M^2 \end{array}\right)$
 is a linear map
\bea
I( ,x) \ :  M^1  \  \rightarrow (Hom(M^2,M^3) \{x \}
\ \
u \mapsto I(u,z)=\sum_{\alpha \in {\bf C}} u_{\alpha} z^{- \alpha-1}
\eea
satisfying the following conditions:
\item[{\bf (I1)}]
For any fixed $u \in M^1$, $v \in M^2, \alpha \in C, u_{\alpha +n}v=0$
for $n$ sufficiently large;
\item [{\bf (I2)}]
$I(L(-1)u,z)v=\frac{d}{dz}I(u,z)v$ for $u \in M^1$, $v \in M^2$;
\item[{\bf (I3)}]  The {\em Jacobi identity} holds, i.e.
\bea
&&z_{0}^{-1}\delta\left(\frac{z_{1}-z_{2}}{z_{0}}\right)
Y_{M^3} (a,z_{1})I (u,z_{2})v
-z_{0}^{-1}\delta\left(\frac{-z_{2}+z_{1}}{z_{0}}\right)
I(u,z_{2})Y_{M^2}(a,z_{1})v=
\nonumber \\
&&z_{2}^{-1}\delta\left(\frac{z_{1}-z_{0}}{z_{2}}\right)
I(Y(a,z_{0})u,z_{2})v
\eea
for any $a \in V$, $u \in M^1$, $v \in M^2$.
\end{definition}

Denote by   $I \left(\begin{array}{c}M^3 \\ M^1, \ M^2 \end{array}\right)$
the vector space of all intertwining operators of this type. Dimension of
this space we call {\em fusion rules} of corresponding type.

By using the definition of fusion rules we can define (under some
circumstances) so called {\em fusion algebra} [DL], [Z]. It's very
interesting that in many examples like WZNW models ([DL], [FZ], [Li1]) 
minimal models ([W], [KWn]) the fusion rules are the same as in CFT 
definition of 'fusion' and also in geometrical definition .

\section{Tensor product}

Let $V_1$ and $V_2$ be two vertex operator algebras.
Following [B] and [FHL] we intoduce the 
VOA structure on vector space $V$=$V_1 \otimes V_2$.

Firstly, $V_1 \otimes V_2$ is on natural way ${\bf Z}_{\geq 0}$-graded
vector space
with natural graduation.
We define the field associated to $a \otimes b \in V$ as follows:
$Y(a \otimes b,z)=Y(a,z) \otimes Y(b,z)$.
If $\omega_1$ and $\omega_2$ are the Virasoro elements of $V_1$
and $V_2$ it's easy to see
that vector $\omega_1 \otimes {\bf 1} + {\bf 1} \otimes \omega_2$
is a Virasoro element for $V$.

The axioms [V1], [V2] and [V3] are satisfied so (with this definition)
$V$ is a VOA.

Now the following simple lemma holds

\begin{lemma}

Let $V_1 \, V_2$ and $V_3$ be a three VOAs. Then :
\item[(1)]
$V_1 \otimes (V_2 \otimes V_3) \cong (V_1 \otimes V_2) \otimes V_3$
\item[(2)]
$V_1 \otimes V_2 \cong V_2 \otimes V_1$

\end{lemma}

Also we need this fact.

\begin{proposition} \label{ppp} [FHL]
Let $M_1$ ($M_2$) be irreducible $V_1$ ($V_2$)--module with $M(0)$ property.
Then $M_1 \otimes M_2$ is irreducible $V_1 \otimes V_2$--module.

\end{proposition}

\section{ The associative algebra $A(V)$ }

\subsection{ Definition and properties }

Let $V$ be a VOA. For any homogeneous
element $a\in V$ and for any $b\in V$, following [Z] we define
\bea
a*b={\rm Res}_{z}\frac{(1+z)^{{\rm wt a}}}{z}Y(a,z)b.
\eea
\bea
a \circ b ={\rm Res}_{z}\frac{(1+z)^{{\rm wt a}}}{z^{2}}Y(a,z)b.
\eea

Then extend this products  bilinearly to the whole space $V$. Let $O(V)$
be a subspace of $V$ linearly spanned with set $\{a \circ b \ : a,b \in V \}$.

Set $A(V)=V/O(V)$. The multiplication $*$ induces the multiplication on
the $A(V)$ and $A(V)$ becomes an associative algebra. The image of $\bf 1$
in $A(V)$ becomes the identity element till the image of $\omega$ is in center of
$A(V)$ (see [Z]).

\subsection{Structure of $A(V_1 \otimes V_2)$}

As we saw in previous section to every VOA we associate the associative
algebra $A(V)$. Here we are going to prove that associative algebras
$A(V_1 \otimes V_2)$ and $A(V_1) \otimes A(V_2)$ are isomorphic.

First we prove some auxiliary results:

\begin{proposition} \label{m}
Let $a,b \in V$ then
$$
a {\circ}_{m}b:={\rm Res}_{z}\frac{(1+z)^{{\rm wt a}}}{z^{m}}Y(a,z)b \in O(V),
$$
for every $m \geq 2$.
\end{proposition}
{\em Proof.}
Put $a'=L(-1)^{m-1}a$. Then it's easy to see that $a' \circ b$=C$(a {\circ}
_{m} b)$ for some constant $C$.

\begin{lemma} \label{ten}

Let $a, a_1  \in V_1$ and $b \in V_2$ then
$$ (a \otimes {\bf 1}) \circ (a_1 \otimes b)=(a \circ a_1) \otimes b$$
\end{lemma}
{\em Proof.}
From the fact ${\bf 1}_{n}a={\delta}_{-1,n}a$ and from the definition of
bilinear mapping ${\circ}$ it follows:
$$(a \otimes {\bf 1}) \circ (a_1 \otimes b)= \sum_{m \in Z}(a \circ_{m} a_1
) \otimes {\bf 1}_{-m+1} b=(a \circ a_1) \otimes b.$$

\begin{lemma} \label{kvoc}
\bea \label {veza}
O(V_1 \otimes V_2) =  O(V_1) \otimes V_2 + V_1 \otimes O(V_2).
\eea
\end{lemma}

{\em Proof.}
From Lemma \ref {ten} follows that $O(V_1) \otimes V_2 + V_1 \otimes O(V_2)
\subseteq O(V_1 \otimes O_2)$. The reverse direction follows from the
definition of $\circ_{m}$ and Proposition {\ref m}.

\begin{theorem} \label {izo}
We define the mapping
\bea \label {formula}
F \ : \ A(V_1 \otimes V_2) \rightarrow A(V_1) \otimes A(V_2), \nn
\eea
\bea \label {pres}
[a \otimes b ] \mapsto [a ] \otimes [b ].
\eea
The mapping $F$ is isomorphism of associative algebras
with identity.
\end{theorem}
{\em Proof.}
Because of 
\bea \label {veza1}
&&A(V_1) \otimes A(V_2)=(V_1 \otimes V_2) /  (O(V_1) \otimes V_2 + V_1 \otimes O(V_2))
\eea
and (\ref {veza}) it follows that this mapping is well-defined.
This map is also surjective. \\
We want to prove injectivity.
Let $m= \sum_{i} [a_i \otimes b_i]$ and $F(m)=0$,
then $ \sum_i a_i \otimes b_i \in O(V_1) \otimes V_2 + V_1 \otimes O(V_2),$
and from (\ref {veza1}) $ \sum_i a_i \otimes b_i \in O(V_1 \otimes V_2)$
follows that $\sum_i [a_i \otimes b_i]=0 ,$ so we prove injectivity.



 
We denote the multiplication in
$A(V_1) \otimes A(V_2)$ with the same simbol $*$.
Now, 
\bea
&&F([a \otimes b]*[a' \otimes b'])=Res_{z}[Y(a \otimes b,z)
\frac{(1+z)^{\rm wt a+wt b}}{z} a' \otimes b'] \nn
&&=[Res_{z}(Y(a,z)\frac{(1+z)^{\rm wt a}}{z}a' \otimes Y(b,z)(1+z)^{\rm wt b}b')] \nn
&&=\sum_{i \in {\bf Z}}[ a {\circ}_{-i} a' \otimes b \circ_{i+2} b'] \nn
&&=(\mbox {by Proposition} {\ref m}) ([a]*[a']) \otimes ([b]*[b']) \nn
&&= ([a]\otimes [b])*([a'] \otimes [b']) \nn
&&=(\mbox{ by Lemma} \ref {kvoc}) F([a \otimes b])*F([a' \otimes b']).
\eea
so $F$ is homomorphism.

\section{Modules and fusion rules for $V_1 \otimes V_2$}

Among all VOAs the most interesting are 'rational VOAs'.

\begin{definition} \label{dd} [Z],[FZ]
VOA is $rational$ if it have only finitely many irreducible representations
and if every finitely generated module for $V$ is completely reducible.
For the VOA which have only finitely many irreducible modules we said that
is $quasirational$.

\end{definition}




\subsection{ Zhu's theorem }

Let $M=\oplus _{n \in {\bf Z} _+ }  M_{n} $ be a $V$--module. For a homogeneous
element $a \in V$ we define $o(a)=a(\deg a-1)$. From the definition of $M$
follows that operator $o(a)$ preserves the gradation of $M$.
\begin{theorem} \label{fz} [FZ]
\item [(a)]
On $\mbox {End}(M_0)$ we have
\bea
&&o(a)o(b) = o(a*b) \nonumber \\
&&o(x) =0 \nonumber
\eea
for every $a,b \in V$, $x \in O(V)$.  The top level $M_0$ is an
$A(V)$--module.
\item[(b)]  Let $U$ be an $A(V)$--module, there exists $V$--module
$M$ such that $A(V)$--module $M_0$ and $U$ are isomorphic.

Thus, we have one-to-one correspondence between irreducible $V$--modules
and irreducible $A(V)$--modules.
\end{theorem}

\begin {remark} \label{vazno}
There is a explicit way how to every $A(V)$--module we can associate
$V$--module.
Morever, if $L(0)$ acts semisimple then for semisimple $A(V)$--module
$M$ corresponding $V$--module is semisimple ([Li],[Z],[DLM]).
\end{remark}


\subsection{Classification theorem}

Put $V=V_1 \otimes V_2$.
In this paragraph we give a the proof of Frenkel-Huang-Lepowsky theorem
of classification of all irreducible $V$--modules with some weaker
condition by using different approach ($A(V)$--theory).

\begin{remark}
We will (from now till the end) assume that $L(0)$ acts semisimply
on every $V$--module.
\end{remark}

Now we need prove technical lemma about representations of
tensor product of two associative algebras with identity.
The proof is from [W].

\begin{lemma} \label{main}

Let $A$ and $B$ be two associative algebras with identity and $M$ a
finite-dimensional irreducible $A \otimes B$--module.
Then $M \cong {M_1 \otimes M_2}$ where $M_1$
($M_2$) is irreducible module for $A$ ($B$).

\end{lemma}

{\em Proof.}
Let $M_1$ be an irreducible submodul of $V$ (for $A$). Then we are
looking the sum $\sum_{b \in B} b {M_1}$ ($b$ means $1 \otimes b$).
This sum is a whole space $M$ because $M$ is irreducible $A \otimes B$
--module. Now $M$ is (as $A$--module) completely reducible and every
irreducible submodule is isomorphic to $M_1$ (this follows from Schur's Lemma).
Set $M_2=Hom_A(M_1,M)$. Then $M_2$ is on natural way $B$-module.
The module $M'=M_{1} \otimes M_{2}$ is $A \otimes B$-- module . 
Now we have $A \otimes B$-homomorphism from  $M'$ 
to a submodule of $M$ such that $$ a \otimes f \mapsto f(a). $$
This homomorphism is bijection so it's isomorphic to $M$.

\begin{remark}
From Burnside's theorem follows that every $A \otimes B$--module
of the form $M_1 \otimes M_2$ is irreducible if $M_1$ and $M_2$
are finite-dimensional. If $M_1$ or $M_2$ is infinite-dimensional
Burnside's theorem doesn't implies the irreducibility.
\end{remark}

The following corollary is a modification of [FHL].
\begin{corollary} \label{co} 
\item [(1.)]
If $V_1$ and $V_2$ are two VOA then every 
irreducible $V$--module M with $M(0)$ property is isomorphic
to $M_1 \otimes M_2$ where $M_1$ and $M_2$ are irreducible $V_1$
($V_2$)--modules.
\item [(2.)]
If $A(V_1)$ and $A(V_2)$ are finite-dimensional 
then every irreducible $V$--module
is of the form $M_1 \otimes M_2$ where $M_1$ and $M_2$ are
irreducible. 
\end{corollary}

{\em Proof.}
Follows directly from Theorem \ref {izo}, Theorem \ref {fz}, 
Lemma \ref {main} and Lemma \ref {ppp}.

\begin{remark}
In some cases of VOA (WZNW models [FZ],[KWn], minimal models [W]) associative algebra
$A(V)$ is a finite dimensional (and rational).
\end{remark}

From the definition of rationality follows:

\begin{lemma} \label {opis}
If $V$ is a rational VOA then $A(V)$ is finite-dimensional semisimple
associative algebra with identity.
\end{lemma}

This fact is also proved in [DLM] with weaker definition of rationality. 


Now we are going to discuss rationality of $V$
if $V_1$ and $V_2$ are rational.
From Remark \ref {vazno} we must analyze when finitely generated 
$A(V)$--module is semisimple. This is equivalent 
to the fact that every $A(V)$--module is semisimple.  
Similar problem was discussed in [DLM1].

\begin{theorem}
If $V_1$ and $V_2$ are two rational VOAs 
then $V=V_1 \otimes V_2$ is also rational. 
If $V_1$ is rational and $A(V_2)$ is finite dimensional
then $V=V_1 \otimes V_2$ is quasirational.
\end{theorem}

{\em Proof.}
Second fact follows directly from the Corrolary \ref {co}.
If $V_1$ and $V_2$ are rational $V$ is from the second fact 
quasirational VOA.
But tensor product of two finite dimensional semisimple algebras 
is also semisimple algebra so every finitely generated
module is semisimple ([P],[J]) and $V$ is rational.

\subsection{Fusion rules}

Let $V_1 \otimes V_2$ be VOA and $M^{i}=M_{1}^{i} \otimes M_{2}^{i}$, $i=1,2,3$
its three irreducible modules such that $L_1(0)$, $L_2(0)$ act semisimply 
on $M_{1}^{i}$, $M_{2}^{i}$ i=1,2,3. 

First, we need the following fact

\begin{proposition} \label {jac}
[Li],[DL]
Let $I(\ ,z)$ be an intertwining operator of type
$\left(\begin{array}{c}M^3 \\ M^1, \ M^2 \end{array}\right)$
then the Jacobi identity is equivalent to following conditions :
For every $a \in V$, $v \in M_1$, $w \in M_2$ exists $r,k \in {\bf N}$ such
that
('associativity')
$$ (z_0 + z_2)^{r}Y(a,z_0)I(v,z_2)w=(z_0 + z_2)^{r}I(Y(a,z_0)v,z_2)w $$
and ('commutativity')
$$ (z_0 - z_1)^{k}Y(a,z_0)I(v,z_2)w=(z_0 - z_1)^{k}I(v,z_2)Y(a,z_0)w $$
are satisfied.
\end{proposition}

Now using this description of 'Jacobi identity' for intertwinig operators
we have the following theorem.

\begin{proposition} \label{manje}

Let $V$ and $M^i$, $i=1,2,3$ be as above then
$$dim \ I \left(\begin{array}{c}M^3 \\ M^1, \ M^2 \end{array}\right) \geq
 dim \ I \left(\begin{array}{c}M_1^3 \\ M_1^1, \ M_1^2 \end{array}\right)
 dim \ I \left(\begin{array}{c}M_2^3 \\ M_2^1, \ M_2^2 \end{array}\right).$$
\end{proposition}
{\em Proof.}
If $I_1$ is an intertwining operator of type
$  \left(\begin{array}{c}M_1^3 \\ M_1^1, \ M_1^2 \end{array}\right)$
and $I_2$ of type
$  \left(\begin{array}{c}M_2^3 \\ M_2^1, \ M_2^2 \end{array}\right)$
then we define $I( \ ,z)$ as on the following way:  
\bea \label {deffusion}
I(m_1 \otimes m_2,z)n_1 \otimes n_2 = I_1(m_1,z)n_1 \otimes I_2(m_2,z)n_2,
\eea
for every $m_1 \otimes m_2 \in M^1$ and $n_1 \otimes n_2 \in M^2$.
Now the first two property ([I1] and [I2])
of intertwining
operators are obviously satisfied. To prove Jacobi identity we use
lemma \ref {jac}. Associativity and commutativity we prove on following
way. 

If associativity for $I_1$ is satisfied for some $r_1 \in {\bf N}$,
$a_1 \in V_1$, $v_1 \in M_1$ and $I_2$ for some $r_2 \in {\bf N}$,
$a_2 \in V_2$, $v_2 \in M_2$. Then for $u=a_1 \otimes a_2$
and $v=v_1 \otimes v_2$ we put $r=max \{ r_1,r_2 \}$ and associativity
is satisfied. In similar way we can prove commutativity.
Now it is easy to see that proposition is true.

The equality in Proposition \ref {manje} 
is not satisfied in general but in
some cases it will hold.

\subsection{The structure of bimodules}

Let $M=M_1 \otimes M_2$ be a $V=V_1 \otimes V_2$--module. Then we define
$A(M)=M/O(M)$ where $O(M)=\{ a \circ m : a \in V, m \in M \}$. We define
left and right action of $A(V)$ on $A(M)$ (see [Z]):
$$ a \cdot m=[Res_z \frac {(1+z)^{wta}}{z}Y(a,z)m ] ,$$ 
$$ m \cdot a=[Res_z \frac {(1+z)^{wta-1}}{z}Y(a,z)m] $$
for every $a \in A(V)$, $m \in A(M)$.
With this notation $A(M)$ is a $A(V)$--bimodule.

\begin{theorem} \label{teh}

Let $M_1$ ($M_2$) be irreducible $V_1$ ($V_2$)-- modules.
Let $a \in A(V_1)$, $b \in A(V_2)$ and $m_1 \in A(M_1)$, $m_2 \in A(M_2)$ then
\item[(1)]
$$
(a \otimes b) \cdot (m_1 \otimes m_2)=(a \cdot m_1) \otimes (b \cdot m_2),  
$$
\item[(2)]
$$
(m_1 \otimes m_2) \cdot (a \otimes b)=(m_1 \cdot a) \otimes (b \cdot m_2). 
$$
Moreover relations (1) and (2) define an isomorphism of
$A(M_1 \otimes M_2)$ and  $A(M_1) \otimes A(M_2)$ as
$A(V_1) \otimes A(V_2)$--bimodules.

\end{theorem}

{\em Proof.} 

From the irreducibility of $M_1$ and $M_2$ follows that
$Y_{M_1}( {\bf 1},z)$ acts as a scalar on $M_1$ (same thing for $M_2$).
Then it's not hard to prove:
  
$$O(M_1 \otimes M_2)= O(M_1) \otimes M_2 + M_1 \otimes O(M_2),$$
using the proofs of Lemmas 4.1 and 4.2.

The proof of (1) is analogous as the proof of Theorem \ref {izo}.
We want to prove (2) which with upper relation
would give us the result.

We have
\bea
&&(m_1 \otimes m_2) \cdot (a \otimes b) \nn  
&&=Res_{z}( \frac {(1+z)^{wta+wtb-1}}{z}Y(a,z)m_1 \otimes Y(b,z)m_2) \nn
&&=\sum_{i \in {\bf Z}} (Res_z (1+z)^{wta-1}z^{-i}Y(a,z)m_1 \otimes Res_{z} (1+z)^{wtb}z^{i-2} Y(b,z)m_2) \nn
&&=\sum_{i>0} (Res_z(1+z)^{wta-1}z^{-i}Y(a,z)m_1 \otimes Res_z (1+z)^{wtb}
z^{i-2}Y(b,z)m_2) \ mod \ O(M_1 \otimes M_2) \nn
&&=Res_z(1+z)^{wta-1}z^{-1}Y(a,z)m_1 \otimes Res_z(1+z)^{wtb}z^{-1}Y(b,z)m_2 \nn
&&+ \sum_{i >1}(Res_z (1+z)^{wta-1}(z^{-i}+z^{-i+1}-z^{-i+1})Y(a,z)m_1 \nn 
&&\otimes Res_z(1+z)^{wtb}z^{i-2}Y(b,z)m_2) \ mod \ O(M_1 \otimes M_2) \nn
&&=m_1 \cdot a \otimes Res_z(1+z)^{wtb}z^{-1}Y(b,z)m_2 \nn
&&+ \sum_{i>1} (Res_z ((1+z)^{wta}z^{-i}-(1+z)^{wta-1}z^{-i+1})Y(a,z)m_1 \nn
&&\otimes Res_z(1+z)^{wtb}z^{i-2}Y(b,z)m_2) \ mod \ O(M_1 \otimes M_2) \nn
&&= m_1 \cdot a \otimes Res_z(1+z)^{wtb}z^{-1}Y(b,z)m_2-Res_z(1+z)^{wta-1}z^{-1}Y(a,z)m_1 \nn
&&\otimes Res_z(1+z)^{wtb}Y(b,z)m_2 - \sum_{i>2}(1+z)^{wta-1}z^{-i+1}Y(a,z)m_1 \nn
&&\otimes Res_z(1+z)^{wtb}z^{i-2}Y(b,z)m_2 \ mod \ O(M_1 \otimes M_2) \nn
&&=m_1 \cdot a \otimes \sum_{i>0}Res_z (1+z)^{wtb}z^{i-2}Y(b,z)m_2 \ mod \ O(M_1 \otimes M_2) \nn
&&=m_1 \cdot a \otimes Res_z(1+z)^{wtb} \frac{1}{z(1+z)}Y(b,z)m_2 \ mod \ O(M_1 \otimes M_2) \nn
&&=(m_1 \cdot a) \otimes (m_2 \cdot b) \ mod \ O(M_1 \otimes M_2) \nn
&&=(m_1 \cdot a) \otimes (m_2 \cdot b). \nn
\eea
In the sumation we use the fact that for $n>>0$ 
$$Res_z(1+z)^{wta}z^{n}Y(a,z)m_1=0,$$
so we have the proof.

We need the following formula for caluclating 'fusion rules'.

\begin{theorem} \label{fusion} [FZ]
Let $M_1, M_2$ and $M_3$ be three V-modules such that $L(0)$ acts semisimple. 
There is one to one corespodance between 
\bea
Hom_{A(V)}(A(M_1) \otimes_{A(V)} M_2(0),M_3(0)) \nn
\eea 
and 'fusion rules'of corresponding type ( $ I \left(\begin{array}
{c}M_3 \\ M_1, M_2 \end{array} \right)$).
\end{theorem}

\begin{lemma} \label{lema}
Let $V$ be a rational VOA $M_1$, $M_2$ and $M_3$ its irreducible 
such that every 'fusion rules' is finite-dimensional then vector space  
\bea
M= A(M_1) \otimes_{A(V)} M_2(0)
\eea
is also finite-dimensional.
\end{lemma}

{\em Proof.}
From the rationality of $V$ follows that $M$ is semisimple and
then can be written:
\bea \label {bez}
M= \oplus_{i \in I} M_i,
\eea   
where $M_i$ are isotypical components  $A(V)$--module $M$. 
But $A(V)$ has only finite number of inequivalent irreducible modules
so set $I$ is finite. If some $M_i$ is infinite-dimensional, and $N$
therm from $M_i$ then $Hom_{A(V)}(M,N)$ would be infinite-dimensional 
which contradict with finitness of 'fusion rules'. 

Now we need following fact from tensor algebra.

\begin{lemma} \label{lemica} 
Let $A$, $B$ be two associative algebras with indentity over ${\bf C}$ and
$M$ and $N$ finite-dimensional right $A$, ($B$)--modules 
and $M'$, ($N'$) finite dimensional left $A$, ($B$)--modules then:
\bea
&&dim  (M \otimes_A M') \otimes \ (N \otimes_B  N'))  \nn
&&=dim (M \otimes M') \otimes_{A \otimes B} (N \otimes N').
\eea
\end{lemma}
Using this Lemma we have:

\begin {theorem}
(Same notation like in Proposition \ref {manje} If $V_1$ and $V_2$ are
rational VOAs with finite 'fusion rules' then in Proposition \ref {manje} 
stands equality.
\end {theorem}

{\em Proof.}

From [KWn] we have:
\bea
&&dim Hom_{A(V)}(A(M_1) \otimes_{A(V)} M_2(0), M_3(0)) \nn
&&=dim M_3(0)^{*} \otimes_{A(V)}A(M_1) \otimes_{A(V)}M_2(0).
\eea
Now if we put 
$A(V)=A(V_1) \otimes A(V_2)$, $M_2(0)=M_1^1(0) \otimes M_2^2(0)$
, $M_3(0)=M_3^1(0) \otimes M_3^2(0)$ and $A(M_1)=A(M_1^1 \otimes M_1^2)$
we have from Lemma \ref {lemica}, \ref {lema}  and Theorem \ref {teh} 
\bea
&&dim \ M_3(0)^* \otimes_{A(V)} A(M_1) \otimes_{A(V)} M_2(0) \nn
&&=dim M_3^1(0)^* \otimes_{A(V_1)} A(M_1^1) \otimes_{A(V_1)} M_2^1(0) \nn
&&dim M_3^2(0)^* \otimes_{A(V_2)} A(M_1^2) \otimes_{A(V_2)} M_2^2(0).
\eea
So we have the proof.

\begin{corollary}
If we have two rational VOAs with Zhu conditions and with 'fusion
rules' 0 or 1 then tensor product have also 'fusion rules' 0 or 1.
\end{corollary}

Department of Mathematics, University of Zagreb,
Bijeni\v{c}ka 30, 41000 Zagreb, Croatia

e-mail :  milas@cromath.math.hr

\end{document}